\documentclass[11pt]{article}
\catcode`\@=11 \catcode`\@=12

\usepackage{amsmath}
\usepackage{amssymb}
\usepackage{amsthm}
\usepackage{oldgerm}

\usepackage[pdftex,colorlinks,urlcolor=blue,pdfstartview=FitH]{hyperref}

\newtheorem{df}{Definition}[section]
\newtheorem{lm}[df]{Lemma}

\newtheorem{Th}{Theorem}

\newtheorem{hyp}{Hypothesis}

\newcommand{\Rm}{\mathbb{R}}

\renewcommand{\SS}{\ensuremath{\mathcal{SS}}}
\newcommand{\SSC}{\ensuremath{\mathcal{SS}_C}}

\newcommand{\Nm}{\ensuremath{\mathbb{N}}}

\newcommand{\mA}{\ensuremath{\mathcal{A}}}

\def\proof {\noindent{\sc{Proof. }}}
\def\qed {\mbox{}\hfill {\small \fbox{}} \\}
\def\lto{\longrightarrow}
\def\lmto{\longmapsto}

\def\leq{\leqslant}
\def\geq{\geqslant}

\newcommand{\e}{\epsilon}

\newcommand{\R}[1]{\mathbb{R}^{#1}}

\newcommand{\ttt}{T_c^-}

\newcommand{\TTT}{T_c^+}

\newcommand{\A}{\mathcal{A}}

\newcommand{\aaaa}{\overline{\mathcal{A}}}

\renewcommand{\d}{{\rm d}}

\author{P. BERNARD, M. ZAVIDOVIQUE}
\title{Regularization of subsolutions   in discrete weak KAM theory}

\begin{document}

\maketitle

\begin{abstract}
We expose different methods of regularizations of subsolutions
in the context of discrete weak KAM theory.
They allow to prove the existence and the density of $C^{1,1}$
subsolutions. Moreover, these subsolutions can be made strict
and smooth outside of the Aubry set.
\end{abstract}

\section{Introduction}

We consider a smooth connected Riemannian manifold $M$
 endowed with the distance $\d(.,.)$ coming from the Riemannian metric.
 Fixing a cost function $c:M \times M\to \R{}$
we study the functions 
$u: M\to \R{}$ which satisfy
$$\forall (x,y)\in M \times M,\ \    u(y)-u(x)\leqslant c(x,y),$$
we call them subsolutions, 
by analogy with those appearing in Weak KAM theory
(see \cite{FaSi,Be1} for example). We will denote by $\SS$ the set of subsolutions, and by $\SSC=\SS\cap C^0(M,\Rm)$ the set of continuous subsolutions.
These subsolutions are one of the important objects in discrete (in time) weak KAM theory.
Some other aspects of this discrete theory have been discussed  in 
\cite{Gom,Be,Za}. This theory is also closely related to the time--periodic
weak KAM theory, discussed for example in \cite{cis,Be111} and many other papers.
In many aspects, these various settings (discrete, time--periodic, autonomous)
are similar, but differences appear for some specific questions.
For example, the convergence of the Lax--Oleinik semi--group holds only in the 
autonomous setting, see
\cite{fathi1998convergence,fathi2000failure,bernard2002connecting,bernard2005convergence},
the Hamilton--Jacobi equation does not have such a nice form
in the discrete setting as in the autonomous setting, see \cite{Gom}. Some other specific
aspects of the discrete case are discussed in \cite{Za}.
Concerning the regularity of subsolutions, 
the existence of $C^1$ subsolutions was obtained in \cite{Za} in the discrete setting 
by an adaptation of the original proof of Fathi and Siconolfi \cite{FaSi}.
On the other hand, the proof of the existence of $C^{1,1}$ subsolutions
given in  \cite{Be1} for the autonomous setting 
does not extend to the discrete setting.
The existence of $C^{1,1}$ subsolutions was however obtained in \cite{Za2} by a different 
method. Our goal here  is to extend and simplify the results of this paper.

Defining, as usual, the   discrete Lax--Oleinik operators
$$
\ttt u(x)=\inf_{y\in M}u(y)+c(y,x),\quad \quad 
\TTT u(x)=\sup_{y\in M}u(y)-c(x,y),
$$
we see that a  function $u$ is a subsolution if and only if one of the equivalent relations is verified:
$$
u\leqslant \ttt u \ \ \ \ \mathrm{or}\ \ \ \ \TTT u \leqslant u.
$$
Note as a consequence that the functions $ \ttt u $ and $\TTT u$ are themselves subsolutions
whenever $u$ is a subsolution.
We will use the following hypothesis on $c$.
More concrete hypotheses implying this one are given below.

\begin{hyp}\label{h1}\rm
For each subsolution $u$, the functions $T_c^-u$ 
and $-T_c^+u$ are locally semiconcave\footnote{Throughout the paper, we call  semiconcave what is sometimes called semiconcave 
with a linear modulus.}. 
\end{hyp} 
Subsolutions do not necessarily exist, and, when they exist, they are not necessarily
continuous (the continuity of subsolutions is discussed in \cite{Za}).
 Under Hypothesis \ref{h1}, the existence
of a continuous subsolution is implied by the existence of a (possibly discontinuous)
subsolution
$u$, just consider   the subsolution $T_c^-u$, which is locally semiconcave
hence locally Lipschitz. See also Lemma \ref{cf} below.

\begin{Th}\label{C}
If Hypothesis \ref{h1} holds, then the set of locally  $C^{1,1}$ subsolutions is dense in the set of
continuous subsolutions for the strong topology.
\end{Th}

We recall that the strong (or Whitney) topology on $C^0(M,\R{})$ is induced by the
basis of open sets:
$$O_{\e,f} = \{ g\in C^0(M,\R{}),\ \ \forall x\in M, \ \ |f(x)-g(x)|< \e (x)\}$$
where $f\in C^0(M,\R{})$ and $\e$ is a continuous positive valued function on $M$. 
For further precisions on this topology, see \cite[Chapter 2]{Hir}.
The existence of $C^{1,1}$ subsolutions was proved in \cite{Za2},
but the density  is new.
In \cite{Za2}, the existence of $C^{1,1}$ subsolutions
is deduced from the following result of Ilmanen
(see  \cite{Il, Carda, fz,  Beil}):

\begin{Th}\label{Il}
Let $f$ and $g$ be locally semiconcave functions
on $M$  such that $f+g\geq 0$.
Then there exists a locally $C^{1,1}$ function $u$ such that 
$-g\leq u\leq f$.
\end{Th}

We will  offer a direct proof of Theorem \ref{C}, which 
is inspired from the proof of Ilmanen's Lemma given in \cite{Beil}.
Note that Theorem \ref{C} implies Theorem \ref{Il}.
This follows immediately from the equivalence, for a given function $u$,
between the two following properties:
 \begin{itemize}
\item
the function $g+u$ is bounded from below and  $-g\leq u-\inf (g+u)\leq f$;
\item
the function $u$ is a subsolution for the cost $c(x,y)=g(x)+f(y)$.
\end{itemize}

We need to introduce more definitions before we  state our other  results.
The subsolution $u$ is called free at $x$ if
$$T_c^+u(x)<u(x)<T_c^-u(x).
$$
We define the  set $\mA_u$
as 
 $$
\mA_u:=\{x\in M,\quad T_c^+u(x)=u(x)=T_c^-u(x)\}\subset M,
$$
and the Aubry set $\mA$ as
$$
\mA:=\bigcap _{u\in \SS} \mA_u\subset M
$$
where the intersection is taken on all  subsolutions.
Under hypothesis \ref{h1},
the sets $\mA_u$ are closed, since they are defined
by the equality $T_c^+u=T_c^-u$.
The set $\mA$ is then also closed. Moreover,
 it makes no difference to restrict the intersection to continuous subsolutions in the definition of $\mA$,
by Lemma \ref{cf} below. 
We say that the subsolution $u$ is strict at $(x,y)$ if
$$
u(y)-u(x)<c(x,y).
$$
Obviously, the subsolution $u$ is strict at $(x,y)$ and 
at $(y,x)$ for each $y$ if it is free at $x$.
We define the set 
$$
\overline \mA_u:=\{(x,y)\in M^2: u(y)-u(x)=c(x,y)\}. 
$$
We also define 
$$
\overline \mA:= \bigcap_{u\in \SS} \overline \mA_u
$$
where the intersection is taken on all subsolutions.
Equivalently, if Hypothesis \ref{h1} holds, the intersection
can be taken on continuous subsolutions, by Lemma \ref{cf}. This yields that $\overline \mA$ is also closed.

\begin{Th}\label{uv}
Assume that $c$ satisfies Hypothesis \ref{h1}.
Given  a subsolution $u$, there exists a subsolution $v$
such that 
\begin{itemize}
\item $v=u$ on $\mA_u$,
\item $v$ is smooth and free on the complement of $\mA_u$,
\item $v$ is locally $C^{1,1}$,
\item $v$ is strict at each pair $(x,y)$ where $u$ is strict.
\end{itemize}
\end{Th}

We can then obtain a subsolution which is  as smooth, free, and strict as possible:

\begin{Th}\label{smoothfree}
If $c$ satisfies hypothesis \ref{h1} and admits  a 
subsolution, then there exists a locally
$C^{1,1}$ subsolution which is free and smooth in the complement of $\A$, and strict on the complement of $\overline \mA$.
\end{Th}

Observe as a consequence that 
the projections of $\overline \mA$ on both the first and the second factor are contained in $\mA$, 
(and, under the additional hypothesis \ref{h2}, each of these projections is equal to $\mA$, see below). 
Strict $C^{1,1}$ subsolutions were obtained in \cite{Za2} under 
an additional twist assumption. We will use a simple trick of 
\cite{Be1} to obtain directly the general result from Theorem \ref{C}.
That the subsolutions can be made smooth outside of $\mA$
is well-known. It will certainly not be a surprise to specialists
that this can be done without destroying the global
$C^{1,1}$ regularity, although we do not know any reference for this 
statement. We prove it using a regularization  
procedure due to De Rham \cite{derh}. This proof also applies
to the ``classical'' (as opposed to discrete) weak KAM theory.

The abstract Hypothesis \ref{h1} holds in 
a more concrete setting, introduced in \cite{Za}:
\begin{hyp}\label{h2}\rm
The function $c$ satisfies the following properties:
\begin{itemize}
 \item\label{unif} \textbf{uniform super-linearity}: for every $k\geqslant 0$,
   there exists $C(k)\in \R{}$ such that 
   $$\forall (x,y)\in M \times M,\quad 
   c(x,y)\geqslant k\d(x,y)-C(k);$$
\item \label{unifb} \textbf{uniform boundedness}: for every $R\in \R{}$, there
  exists $A(R)\in \R{}$ such that 
 $$\forall (x,y)\in M\times M, \quad  \d(x,y)\leqslant R \Rightarrow
  c(x,y)\leqslant A(R);$$
\item \textbf{local semiconcavity}:
for each point $(x_0,y_0)$ there is a domain of chart 
containing  $(x_0,y_0)$ and a smooth function $f(x,y)$
such that $c-f$ is concave in the chart.
(This holds for example if $c$ is $C^2$ or locally $C^{1,1}$).
\end{itemize}
\end{hyp}

This hypothesis has two important consequences, as was proved in  \cite{Za}. First, it implies
Hypothesis \ref{h1}. Second, it implies that the extrema
in the definitions of 
$T_c^{\pm}u(x)$ are reached for each continuous subsolution
$u$ and each $x\in M$.
This in turn implies that the the projection of $\overline\mA$
on the first, as well as on the second, factor are equal to
$\mA$, which corresponds to the projected Aubry set introduced in \cite{Za}:
\begin{lm}
Assume that $c$ satisfies Hypothesis \ref{h2}.
Given $x\in \A$,  there exist $y$ and $z$ such that $(x,z)$ and $(y,x)$ are in $\aaaa$.
\end{lm}
\proof
Let  $w$ be a continuous  subsolution which is strict outside of
$\overline \mA$ (such a solution exists by Theorem 
\ref{smoothfree}). Let $y$ be such that $\ttt w(x)=w(y)+c(y,x)$.
Since $x\in \A$ we obtain that $w(x)-w(y)=c(y,x)$. Hence $(y,x)\in \aaaa_{w}=\aaaa$. The existence of $z$ is proved in the same way, using $\TTT$.
\qed

Finally, let us mention one last setting in which Hypothesis \ref{h1} holds :
\begin{hyp}\label{h3}\rm
The function $c$ is locally bi--semiconcave:

 for all $(x,y)\in M\times M$ we can find the following:
\begin{itemize}
\item  neighborhoods  $U$ and $V$ of respectively $x$ and $y$, 
\item diffeomorphisms $\varphi_1$ and $\varphi_2$ from $B_n$ to respectively $U$ and $V$ ($B_n$ is the unit ball in $\R{n}$),
\item smooth functions $f$ and $g$ from $B_n$ to $\R{}$,
\end{itemize}
such that for each $x\in M$, the function $z\mapsto c\big(x,\varphi_2(z)\big)-g(z)$ is concave and
for all $y\in M$, the function $z\mapsto c\big(\varphi_1(z),y\big)-f(z)$ is concave.
\end{hyp}
It is  easy to prove, as in  \cite[Proposition 4.6]{Za2},  that Hypothesis \ref{h3} also implies Hypothesis \ref{h1} (using that an infimum of equi--semiconcave functions is itself semiconcave).

\section{Preliminaries}

We gather here some useful facts obtained from 
elementary manipulations of the Lax--Oleinik operators.
Let us first  list, without proof,  some   properties of the operators $T_c^{\pm}$.
\begin{itemize}
\item Monotony : 
$u\leq v \Rightarrow T_c^{\pm}u \leq  T_c^{\pm}v$.
\item Convexity : Given a sequence $u_n$ of functions
and a sequence $a_n$ of non--negative numbers such that 
$\sum_{n\in \Nm} a_n =1$, and such that the series
$\sum_{n\in \Nm}a_nT_c^-u_n$, 
$\sum_{n\in \Nm}a_nu_n$ and $\sum_{n\in \Nm}a_nT_c^+u_n$
 are converging point--wise,  we have
$$
T_c^- \big(\sum_{n\in \Nm}a_nu_n\big)\geq 
\sum_{n\in \Nm}a_nT_c^-u_n
\quad,\quad
T_c^+ \big(\sum_{n\in \Nm}a_nu_n\big)\leq 
\sum_{n\in \Nm}a_nT_c^+u_n.
$$
The set $\SS$ of subsolutions is convex, and it is closed under point--wise convergence.
A convex combination $\sum_{n\in \Nm}a_nu_n$ of subsolutions, with a point--wise convergent sum,
is a subsolution; it is free at $x$
(resp. strict at $(x,y)$) provided there exists $n$ such that $a_n>0$ and such that $u_n$ is free at $x$ (resp. strict at $(x,y)$).
\item We have the equalities $T_c^+\circ T_c^-\circ T_c^+=T_c^+$ and
$T_c^-\circ T_c^+\circ T_c^-=T_c^-$.
\item
We have the inequalities
$$
T_c^+\circ T_c^-u\leq u\quad, \quad T_c^-\circ T_c^+ u\geq u
$$
for each function $u$.
\item If $u$ is a subsolution, then

\begin{equation}\label{ineq}
T^+_c u\leq T^+_c\circ T^-_c u\leq 
u\leq T^-_c\circ T^+_c u \leq T^-_c u
\end{equation}
\end{itemize}

The following criterion 
for subsolutions is taken from \cite{Za}:

\begin{lm}\label{entre}
 Let $u$ be a subsolution and let us consider  a function $v$ such that
$$u\leqslant v\leqslant \ttt u,$$
then $v$ itself is a  subsolution.
\end{lm}
\proof
The statement follows from the inequalities
$ u\leqslant v\leqslant \ttt u \leqslant \ttt v.$
\qed

Playing with the Lax--Oleinik operators also leads to: 
\begin{lm}\label{cf}
Let $u$ be a  subsolution, then the subsolution
$$
v:=\frac{T_c^+u+T_c^+\circ T_c^-u+T_c^-\circ T_c^+u
+T_c^-u}{4}
$$ 
is free on the complement of $\mA_u$,
equal to $u$ on $\mA_u$, and strict on the complement
of $\overline{\mA}_u$. If Hypothesis \ref{h1}
holds, then $v$ is locally Lipschitz.
\end{lm}

We then have $\mA_v\subset \mA_u$, but this inclusion is not
necessarily an equality.

\proof
To prove that $v$ is free on the complement of 
$\mA_u$,
we consider a point $x$ at which $v$ is not free, and prove that $x\in \mA_u$.
We either have  $T_c^+v(x)=v(x)$
or $T_c^-v(x)=v(x)$.
In the first case, we have
$$
4v(x)=4T_c^+v(x)\leq
 T_c^+\circ T_c^+u(x)+T_c^+\circ T_c^+\circ T_c^-u(x)
+T_c^+\circ T_c^-\circ T_c^+u(x)+T_c^+\circ T_c^-u(x)
$$
hence 
the  inequalities
\begin{align*}
T_c^+\circ T_c^+u(x)\leq T_c^+u(x)
\quad,\quad 
T_c^+\circ T_c^+\circ T_c^-u(x)\leq T_c^+\circ T^-_cu(x)\\
T_c^+\circ T_c^-\circ T_c^+u(x)=
T_c^+u(x)\leq T_c^-\circ T_c^+u(x)
\quad, \quad
T_c^+\circ T_c^-u(x)\leq T_c^-u(x)
\end{align*}
sum to an equality, hence they are equalities.
In view of (\ref{ineq}) the two last equalities imply that 
$T_c^+u(x)=u(x)=T_c^-u(x)$.
The second case is similar.
It then follows from Lemma \ref{strict} below that $v$ is strict
outside of $\overline \mA_u$.
\qed

The following Lemma allows to reduce strictness questions to freedom questions, and ends the proof of Lemma \ref{cf}.

\begin{lm}\label{strict}
Let $u,v$ be  subsolutions, such that  $v$
is free outside of $\mA_u$ and equal to $u$ on
$\mA_u$, then $v$ is strict at each point $(x,y)$
where $u$ is strict.
\end{lm}

\proof
Let $(x,y)$ be a pair at which $v$ is not strict.
Then $v(y)-v(x)=c(x,y)$, hence $T^-_cv(y)=v(y)$
and $T^+_cv(x)=v(x)$. Since $v$ is free outside of 
$\mA_u$, this implies that both $x$ and $y$ belong to
$\mA_u$. Since $u=v$ on $\mA_u$, we conclude that 
$$u(y)-u(x)=v(y)-v(x)=c(x,y)
$$
hence $u$ is not strict at $(x,y)$.
\qed

It will also be useful to quantify the freedom of a
 subsolution $u$ by its leverage function:
        \begin{df}\rm
 The leverage function $\lambda_u:M\lto[0,\infty)$ of the subsolution $u$ is defined by:
    $$\lambda_u(x):=\frac{1}{3}\min \big(\ttt u(x)- u(x),\ u(x)-\TTT u(x)\big).$$
 \end{df}
 
Note that $u$ is free at $x$ if and only if $\lambda_u(x)>0$.

 \begin{lm}\label{leverage}
 Let $u$ be a  subsolution and let  $v$ be another function such  that 
 $|u-v|\leqslant \lambda_u,$
 then $v$ is itself a subsolution. Moreover, if $u$ is free at $x$ then so is $v$,
 and if $u$ is strict at $(x,y)$, then so is $v$.
 \end{lm}
 
 \proof
 By definition, we have
\begin{eqnarray*}3\max \{\lambda_u(x),\   \lambda_u(y)\}\leqslant 
 \max \{u(x)-\TTT u(x),\   \ttt u(y)-u(y)\}
\leqslant &c(x,y)-u(y)+u(x) .
\end{eqnarray*}
We conclude that
$$
0\leq \max \{\lambda_u(x),\   \lambda_u(y)\}\leq c(x,y)-v(y)+v(x)
$$
hence that 
$$T_c^+v(x)+\lambda_u(x)\leq v(x)
\leq T_c^- v(x)-\lambda_u(x)
$$
for each $x$, which implies that $v$ is a subsolution 
which is free at points where $u$ is free.
The last claim follows from Lemma \ref{strict}.
\qed  

%
%

\section{The uniform case on $\R{n}$ and the Jensen transforms}

In this section we work on   $M=\R{n}$.
A function $u:\Rm^n \lto \Rm$ is called
$k$--semiconcave if $u(x)-k\|x\|^2$ is concave.
We make the following more quantitative version of Hypothesis \ref{h1} on the cost $c$:

\vspace{3mm}
 \noindent{\bf Hypothesis 1--K.}
There exists a constant $K$ such that for each subsolution $u$, the functions $T_c^-u$ 
and $-T_c^+u$ are $K$--semiconcave.

\vspace{3mm}

One setting which implies this condition is tho following version of Hypothesis \ref{h3}:

\vspace{3mm}
 \noindent{\bf Hypothesis 3--K.}
There exists a constant $K$ such that the function
$x\lmto c(x,y)$ is $K$--semiconcave for each $y$ and 
the function  $y\lmto c(x,y)$ is $K$--semiconcave for each $x$.

\vspace{3mm}

We will use the Jensen transforms which associate,
to a function
$u: \R{n} \to \R{}$ and a positive real number $t$, the functions
 $$J^{-t}u(x)=\inf_{y\in \R{n}} \Big( u(y)+\frac{1}{t}\|y-x\|^2\Big),
 \qquad J^{+t}u(x)=\sup_{y\in \R{n}} \Big( u(y)-\frac{1}{t}\|y-x\|^2\Big).
$$
These are nothing but the Lax--Oleinik operators associated 
to the costs $c_t(x,y)=\frac{1}{t}\|y-x\|^2$.
\begin{Th}\label{Casunif}
Let $u$ be a  uniformly continuous subsolution. The function 
$J^{-t}\circ J^{+2t}\circ J^{-t} u$ is finite and, 
 for $t$ small enough, it is a $C^{1,1}$ subsolution. Moreover, it converges uniformly to $u$ as $t\to 0$.
More precisely,  if $u$ is a uniformly continuous subsolution 
then for $t,s<K^{-1}$ the functions 
$J^{-t}\circ J^{+(t+s)} \circ J^{-s} u$ and
$J^{+t}\circ J^{-(t+s)} \circ J^{+s} u$ are  $C^{1,1}$ 
subsolutions which converge uniformly to $u$ as $t,s\to 0$.
Moreover, we have 
$$T^+_c\circ T_c^-u\leq J^{-t}\circ J^{+(t+s)} \circ J^{-s} u \leq T^-_cu \quad ,\quad
T^+_cu\leq J^{+t}\circ J^{-(t+s)} \circ J^{+s} u \leq T_c^-\circ T_c^+u.
$$
\end{Th}

Note that the last inequalities imply that 
$J^{-t}\circ J^{+(t+s)} \circ J^{-s}u$ and 
$J^{+t}\circ J^{-(t+s)} \circ J^{+s} u$ are subsolutions,
by Lemma \ref{entre}.
We recall a few properties of the Jensen transforms, most of 
which are  proved in \cite{Beil} or \cite{Amb}. 
Both families of operators $J^-$ and $J^+$ are semi--groups. They are monotonous in the following way:
$$\forall s>t>0, \qquad \inf  u \leqslant J^{-s}u \leqslant J^{-t} u \leqslant u \leqslant J^{+t} u \leqslant J^{+s} u \leqslant \sup u$$
and in the following one:
$$u \leqslant v \Rightarrow \{\forall t \geqslant 0,\quad J^{-t}u \leqslant J^{-t}v \; \;\mathrm{and}\;\; J^{+t}u \leqslant J^{+t}v\}.$$
We call modulus of continuity a continuous function 
$\rho:[0,\infty)\lto [0, \infty)$ such that $\rho(0)=0$.
A function $f$ is said $\rho$--continuous if 
$|f(y)-f(x)|\leq \rho(\|y-x\|)$
for all $x$ and $y$.
Given a modulus of continuity $\rho$, there exists a modulus of continuity $\epsilon$ such that, for each $\rho$--continuous function $u$, the following properties hold:
\begin{itemize}
\item 
the functions $J^{-t}u$ and $J^{+t}u$ are finite-valued
and $\rho$--continuous for each $t\geq 0$,
\item  $J^{-t}u$ is $t^{-1}$--semiconcave and $J^{+t}u$ is $t^{-1}$--semiconvex,
\item 
 $\|J^{-t}u-u\|_{\infty}+\|J^{+t}u-u\|_{\infty}\leq \epsilon(t)$,
\item
 $J^{-t}\circ J^{+t}u\geqslant u$ and  $J^{+t}\circ J^{-t}u\leqslant u$,
\item
 the equality $J^{-t}\circ J^{+t}u= u$ (resp. $J^{+t}\circ J^{-t}u= u$) holds if and only if $u$ is $t^{-1}$--semiconcave (resp. $t^{-1}$--semiconvex),
\item
 if $u$ is semiconvex (resp.  semiconcave) then $J^{-t}\circ J^{+t}u$ (resp. $J^{+t}\circ J^{-t}u$) is $C^{1,1}$ (and finite valued).
\end{itemize}
Using these properties, we now prove Theorem \ref{Casunif}.
Let $u$ be a uniformly continuous subsolution, with modulus $\rho$.
Since the function $u$ is a subsolution, we have $u\leq T^-_cu$
hence $T^-_cu$ is finite--valued.
Our hypothesis is that 
the function $T^-_cu$ is $K$--semiconcave.
For $s<K^{-1}$, we have
$$
u\leqslant J^{-s}\circ J^{+s}u \leqslant J^{-s}\circ J^{+s}(\ttt u )=\ttt u,
$$
where the last inequality follows from the $K$--semiconcavity
of $T^-_cu$ and the properties of $J^{-}\circ J^{+}$ listed above.
We conclude that the function 
$J^{-s}\circ J^{+s}u$ is a $\rho$--continuous, $s^{-1}$--semiconcave
subsolution.
Similarly, if $u$ is $\rho$--continuous and $t<K^{-1}$,
then the function 
$J^{+t}\circ J^{-t}u$ is a $\rho$--continuous, $t^{-1}$--semiconvex
subsolution.
Applying this observation to the function $J^{-s}\circ J^{+s}u$,
we conclude that $J^{+t}\circ J^{-t}\circ J^{-s}\circ J^{+s}u$
is a $\rho$--continuous subsolution.
This subsolution is 
$C^{1,1}$ since $J^{-s}\circ J^{+s}u$ is semiconcave.
We have the inequality
$$
T^+_c\circ T_c^-u=J^{+t}\circ J^{-t}(T^+_c\circ T_c^-u)\leq J^{+t}\circ J^{-t} u
\leq J^{+t}\circ J^{-t} \circ J^{-s}\circ J^{+s}u
\leq J^{+t}\circ J^{-t}(T^-_cu)\leq T^-_cu.
$$
Finally, we have
$$
u\leq J^{-s}\circ J^{+s} u \leq J^{-s} (u+\|J^{+s}u-u\|_{\infty})\leq 
\|J^{+s}u-u\|_{\infty}+\|J^{-s}u-u\|_{\infty}+u\leq u+ \epsilon(s)
$$
and similarly
$
u-\epsilon(t)\leq J^{+t}\circ J^{-t}u\leq u
$
hence 
$$
u-\epsilon(t)\leq
 J^{+t}\circ J^{-t}u\leq J^{+t}\circ J^{-(t+s)} \circ J^{+s} u\leq 
 J^{-s}\circ J^{+s} u\leq u+\epsilon(s),
$$
where $\epsilon$ is the modulus associated to $\rho$
in the list of properties of $J$.
\qed

\section{The general case}\label{sectiongen}
 
 In this section, we come back to the general setting
and prove Theorem \ref{C}.
We derive it from the uniform version using partitions
of unity, as was done in \cite{Beil} 
for Ilmanen's Lemma.
We fix a locally finite atlas $(\phi_i)_{i\in I}$ constituted of smooth maps $\phi_i : B_n\to M$, where $B_n$ is the open unit ball. We assume  that all the images $\phi_i(B_n)$, for $ i\in I$, are relatively compact in $M$. Moreover, we consider a smooth partition of unity $(g_i)_{i\in I}$ subordinated to the 
locally finite open covering $\big(\phi_i(B_n)\big)_{i\in I}$.
Given positive numbers $a_i,b_i,  i\in I$, we define the operators \begin{equation}\label{operatorR}
 \forall x\in M, \ \ \ \      S u(x)=\sum_{i\in I}\ [ J^{-a_i}\circ J^{+ a_i}(g_i u\circ \phi_i) ]\circ \phi_i^{-1}(x),
 \end{equation}
  \begin{equation}\label{operRR}
 \forall x\in M, \ \ \ \check S u(x)=\sum_{i\in I} \    [ J^{+b_i}\circ J^{- b_i}(g_i u\circ \phi_i) ]\circ \phi_i^{-1}(x).
 \end{equation}
The functions in the sums are extended to the whole of $M$  by the value zero outside of the domain $\phi_i(B_n)$. The sums are locally finite hence well-defined.
Theorem \ref{C} follows from: 
\begin{Th}\label{operR}
Assume that the cost $c$ satisfies Hypothesis \ref{h1}. 
 Let $u$ be a continuous subsolution and $\e :M\to {]0,\infty)}$ be a continuous function. For suitably chosen positive constants $(a_i)_{i\in I}$ and $(b_i)_{i\in I}$, the function $\check S \circ S (u)$ is a
locally  $C^{1,1}$ subsolution such that $|u-\check S \circ S u|\leqslant \e$
and
$$
T^+_c\circ T_c^-u\leq \check S \circ S u \leq T^-_cu.
$$
 \end{Th}
 
\proof
Since the image $\phi_i(B_n)$ is relatively compact and since the atlas is locally finite
the set $A_i=\{j\in I,\ \ \phi_j(B_n)\cap \phi_i(B_n)\neq \varnothing\}$ is finite, let us denote by $e_i$ its cardinal. 
Setting

$$
\epsilon_i:=\frac{\min\limits_{j\in A_i}\inf\limits_{x\in B_n} \e \big(\phi_j(x)\big)}
{2\max\limits_{j\in A_i}e_j},
$$
we observe  that  
\begin{equation}\label{petit}
\forall i\in I,\ \sum_{j\in A_i} \e_j\leqslant \frac{1}{2}\inf_{x\in B_n} \e \big(\phi_i(x)\big).
\end{equation}
Let us make the convention to extend all functions which are compactly supported
inside $B_n$, like $(g_iu)\circ \phi_i$  by the value $0$ to the whole of $\Rm^n$.
For each $i$, we choose a positive constant $a_i$ such that
\begin{equation}\label{petitt}
 \big\| (g_i u )\circ \phi_i -J^{-a_i}\circ J^{+a_i}\big((g_i u )\circ \phi_i\big)\big\|_\infty <\e_i.
\end{equation}
Such a constant exists because the function 
$(g_iu)\circ \phi_i$ is uniformly continuous on $\Rm^n$.
Since 
 $\ttt u$  is locally semiconcave,
the function   $(g_i \ttt u )\circ \phi_i$, extended by zero outside of $B_n$,
is semiconcave on $\Rm^n$ (see \cite{Beil}).
We can assume by taking $a_i>0$ small enough
that it is $a^{-1}_i$--semiconcave, so that
$$[g_iu]\circ \phi_i \leqslant  J^{-a_i}\circ J^{+a_i}\big([g_i u]\circ \phi_i\big)\leqslant
 J^{-a_i}\circ J^{+a_i}\big([g_i \ttt u]\circ \phi_i\big)=\left[g_i \ttt u\right]\circ \phi_i
$$
on $\Rm^n$. This implies in particular that 
the function  $J^{-a_i}\circ J^{+a_i}(g_i u\circ \phi_i)$ is supported in $B_n$.
As a consequence, the function 
$\big[J^{-a_i}\circ J^{+a_i}(g_i u\circ \phi_i)\big]\circ \phi_i^{-1}$,
extended by zero outside of $\phi_i(B_n)$,
is locally semiconcave on $M$, hence the function $Su$ is locally semiconcave, being
a locally finite  sum of locally semiconcave functions.
By summation, we get
$$u=\sum_{i\in I} (g_i u)\circ \phi_i\circ \phi_i^{-1} \leqslant S u \leqslant \sum_{i\in I} \left[g_i \ttt u\right]\circ \phi_i\circ \phi_i^{-1}=\ttt u,$$
which, by Lemma \ref{entre}, implies that 
 $S u $ is a  subsolution. We have $|u-S u|< \e/2$, by (\ref{petitt}).

Next, we chose $b_i$ such that 
 $[g_i T_c^+\circ T_c^- u]\circ \phi_i$ is 
$b_i^{-1}$--semiconvex, which implies that
\begin{align*}
[g_i T_c^+\circ T_c^- u]\circ \phi_i&=J^{+b_i}\circ J^{-b_i} \big([g_i T_c^+\circ T_c^- u]\circ \phi_i\big)\\
&\leq 
J^{+b_i}\circ J^{-b_i}\big ([g_iu]\circ \phi_i\big)
\leq J^{+b_i}\circ J^{-b_i}\big ([g_iSu]\circ \phi_i\big)
\leq [g_iSu]\circ \phi_i.
\end{align*}
As above, this implies that 
$J^{+b_i}\circ J^{-b_i}\big ([g_iSu]\circ \phi_i\big)$
is supported on $B_n$. Note that it is also $C^{1,1}$ hence
the function 
$\big(J^{+b_i}\circ J^{-b_i}\big ([g_iSu]\circ \phi_i\big)\big)\circ \phi_i^{-1},
$
extended by zero outside of $\phi_i(B_n)$, is locally $C^{1,1}$ on $M$.
By summation, we obtain that
$$
T_c^+\circ T_c^- u\leq \check S u \leq \check S \circ S u \leq Su\leq T^-_cu,
$$
which implies that $\check S \circ S u$ is a subsolution. This function is 
locally $C^{1,1}$ as a locally finite sum of locally $C^{1,1}$ functions.
Finally, we can assume by possibly reducing $b_i$ that
$$
 \big\| (g_i S u )\circ \phi_i -J^{+b_i}\circ J^{-b_i}\big((g_i S u )\circ \phi_i\big)\big\|_\infty <\e_i,
$$
 which implies that $|\check S\circ Su-Su|\leq \e/2$
hence that $|\check S\circ Su-Su|\leq \e$.
 \qed

\begin{Th}\label{Rfree}
We assume hypothesis \ref{h1}.
Let $\Omega\subset M$ be an open set and let $u$ be a continuous 
subsolution which is free on $\Omega$. Then the subsolution $u$
belongs to the closure, for the strong topology,
of 
the set of $C^{1,1}$ subsolutions which are free on $\Omega$ and equal to $u$ on
$\mA_u$.
\end{Th}

\proof
Let $\epsilon:M\rightarrow {]0, \infty)}$ be a continuous function.
We can chose  $a_i$ and $b_i$ in such a way that 
$\check S\circ Su$ is a subsolution which is equal to $u$
on $\mA_u$, and such that $|\check S\circ Su-u|\leq \e$.
However, $\check S\circ Su$ need not be free on $\Omega$.
To preserve the freedom of $u$, we work with  the modified cost
$$\tilde c(x,y)=c(x,y)-\psi(y),
$$
 where
 $\psi$ is a smooth bounded function such that 
$0\leq \psi \leq \lambda_u$ (the leverage function of $u$),
 with strict inequalities on $\Omega$.
The associated Lax--Oleinik operator is 
$$
T_{\tilde c}^- v(x)=-\psi(x)+T_c^-v(x).
$$
Each subsolution for the cost $\tilde c$
is thus  a subsolution for the cost $c$, and
$\tilde c$ satisfies Hypothesis \ref{h1}.
Moreover, 
the function $u$ is a subsolution for the cost
$\tilde c$.
We apply Theorem \ref{operR} and get a locally $C^{1,1}$ subsolution
$w^-$ for the cost $\tilde c$, which satisfies $|w^--u|\leq \epsilon$ and $w^-=u$ on $\mA_u$.
This function then satisfies
$$
T^-_cw^-=\psi+T^-_{\tilde c} w^-\geq\psi+ w^-
$$
hence it is a subsolution for the cost $c$.
Similarly, by applying Theorem \ref{operR} with the modified cost
$c(x,y)-\psi(x)$, we get a locally $C^{1,1}$ subsolution 
$w^+$ (for the cost $c$) such that 
$$
T^+_cw^+\leq w ^+-\psi,
$$
$|w^+-u|\leq \epsilon$ and $w^+=u$ on $\mA_u$.
We then set $w:= (w^++w^-)/2$
and claim that this locally $C^{1,1}$ subsolution 
is free on $\Omega$.
This follows from the inequalities
\begin{align*}
T_c^-w&\geq \big(T_c^-w^-+T_c^-w^+\big)/2
\geq w+\psi/2,\\
T_c^+w&\leq \big(T_c^+w^-+T_c^+w^+\big)/2
\leq w-\psi/2,
\end{align*}
 since $\psi$ is positive  on $\Omega$.
We also obviously have 
$|w-u|\leq\epsilon$ and $w=u$ on $\mA_u$.
\qed

\section{Proof of Theorem \ref{uv}}

We will build  successively subsolutions $v_1,v_2,v_3$
which are all equal to $u$ on $\mA_u$ and free on the 
complement $\Omega$ of $\mA_u$. By Lemma \ref{strict}, this also
implies that the subsolutions $v_i$ are strict where $u$
is strict.
We take 
$$
v_1=\frac{T_c^+u+T_c^+\circ T_c^-u+T_c^-\circ T_c^+u
+T_c^-u}{4},
$$
which is continuous, equal to $u$ on $\mA_u$ and free on the 
complement  of $\mA_u$ by Lemma \ref{cf}.

We then build $v_2$ by applying Theorem \ref{Rfree} to $v_1$, with $\Omega=M\setminus \mA_u$, and
get a locally $C^{1,1}$ subsolution $v_2$ which is free on $\Omega$ and equal to $u$ on $\mA_u$.

The following mollification  result,  
which will be proved in the Appendix using a procedure due to De Rham, allows to smooth our subsolution on $\Omega$.

\begin{Th}\label{approx}
Let $f$ be a locally $C^{k,1}$ function on $M$
 and let $\e:M\lto [0,\infty)$ be a continuous function.
 Then, there exists a locally  $C^{k,1}$ function $g:M\to \R{}$ 
which is smooth on the open set  $\Omega:=\e^{-1}(0, +\infty)$ and satisfies,
for all $x\in M$,
$$  \ \ |f(x)-g(x)|+\|\d_x f - \d_x g\| +\cdots +\| \d^k_x f-\d^k_x g \|\leqslant \e(x).$$
\end{Th}

More precisely, we apply Theorem \ref{approx}
to the function $f=v_2$, with $k=1$,
and with a function $\e(x)$ such that $\e=0$ on $\mA_u$,
$\e>0$ on $\Omega$ (the complement of $\mA_u$), and $\e \leq \lambda_{v_2}$
(the leverage function of $v_2$).
We get a $C^{1,1}$ function $v_3$, which is smooth on
$\Omega$ and is equal to $u$ on $\mA_u$.
Since $|v_3-v_2|\leq\lambda_{v_2}$, Lemma \ref{leverage} implies that $v_3$ is a subsolution which is free on $\Omega$. Lemma \ref{strict} then implies that $v_3$ is strict
where $u$ is strict.
\qed

\section{Proof of Theorem \ref{smoothfree}}
It is enough to prove the existence of 
a subsolution $u$
which is free on the complement of $\mA$
and strict on the complement of $\overline \mA$.
Theorem \ref{uv} then implies the existence of a 
locally $C^{1,1}$ solution
$v$ which is free and smooth on the complement of $\mA$,
and which is strict on the complement of $\overline \mA$.
We start with:

\begin{lm}\label{optimal}
If $c$ satisfies Hypothesis \ref{h1} and  admits a  subsolution,
then there exists a continuous subsolution $w_1$
which is free on the complement of $\mA$.
\end{lm}

\proof
Let us consider a point  $x\notin \A$. By definition,
 there exists a  subsolution $v_x$ 
such that $x\notin  \mA_{v_x}$, hence, by Lemma \ref{cf}, there exists
a continuous subsolution $u_x\in \SSC$ which is free at $x$.
By continuity of $u_x$,  $\ttt u_x$ and   $\TTT u_x$
we may consider a positive number $\e_x$ and an open neighborhood of $x$, $O_x$, 
on which the following holds:
$$\forall y\in O_x, \ \ \ \ \ttt u_x (y) -\e_x> u_x (y)>\TTT u_x(y)+\e_x.$$
The set  $M\setminus \A$ satisfies the Lindel\" of property
(it is a separable metric space). We can  thus extract a countable covering
 $O_n$, $n\in \mathbb{N}$ of the covering $O_x$,  $x\in M\setminus \A$. 
Denoting by $u_n$ and $\e_n$ the continuous subsolution and positive 
real number associated to $O_n$, we consider  a convex combination
$$w_1= \sum_{n\in \mathbb{N}} a_n u_n,
$$
where $a_n$ is a sequence of positive numbers such that $\sum_{\Nm}a_n=1$
and such that the sum in the definition of $w_1$ is normally convergent on each compact set.
The function $w_1$ is then a continuous subsolution.
For each $x\notin \A$, there exists  $n_0\in \mathbb{N}$ such that $x\in O_{n_0}$, and we have
$$\ttt w_1 (x)= \ttt \Big(\sum_{n\in \mathbb{N}} a_n u_n\Big) (x) 
\geqslant \sum_{n\in \mathbb{N}} a_n    \ttt u_n (x) 
\geq a_{n_0} \e_{n_0}+\sum_{n\in \mathbb{N}}a_n u_n  >w_1(x).
$$ 
A similar computation shows that $\TTT w_1(x)<w_1(x)$.
\qed

\begin{lm}\label{fonctionstricte}
If there exists a continuous subsolution, then there exists a continuous subsolution $w_2$ which is strict at each 
pair $(x,y)$ 
where a strict continuous subsolution exists. Under Hypothesis \ref{h1}, the subsolution $w_2$ is then strict outside 
of $\overline \mA$.
\end{lm}
\proof
Since $M$ is separable, the set $\SSC$ of continuous subsolutions is also separable (for the compact--open topology), and we consider
   a dense subsequence $(u_n)_{n\in \mathbb{N}}$. Set 
\begin{equation}\label{series}
w_2=\sum_{n\in \mathbb{N}} a_n u_n
\end{equation}
where the $a_n$ are positive real numbers such that $\sum a_n=1$ and the sum (\ref{series}) is uniformly convergent on each compact subset. The function $w_2$ is a subsolution since it is a convex combination of subsolutions. If now $(x,y)\in \aaaa_{w_2}$,
summing the inequalities
$$\forall n\in \mathbb{N}, \ \ a_n\big(u_n(y)-u_n(x)\big)\leqslant a_n c(x,y),$$
gives an equality, therefore all inequalities are equalities and 
$$\forall n\in \mathbb{N}, \ \ (x,y)\in \aaaa_{u_n}.$$
By density of  the sequence $u_n$,
we deduce that  $(x,y)\in \aaaa_u$ 
for each continuous solution $u$. 
Under Hypothesis \ref{h1},   $\overline \mA$ 
is exactly the set of pairs at which no continuous subsolution is strict: $\overline\mA = \bigcap_{u\in \SSC} \overline \mA_u$ hence, $(x,y)\in \overline \mA$.
\qed

To finish the proof of Theorem \ref{smoothfree},
we consider the subsolution $u=(w_1+w_2)/2$.
This subsolution is free on the complement of $\mA$
because $w_1$ is, and it is strict on the complement of $\overline \mA$ because $w_2$ is.

\qed

\appendix

\section{Proof of  Theorem \ref{approx}}\label{appen}
We prove Theorem \ref{approx} using a regularization procedure  due to De Rham,
see \cite{derh}.
The idea of De Rham is to construct an action $\mathfrak{t}$
of $\Rm^n$ on $\Rm^n$ by smooth diffeomorphisms supported on the unit sphere 
$B_n$, in such a way that the induced action on $B_n$
is conjugated to the standard action of $\Rm^n$
on itself by translations.
More precisely, there exists a diffeomorphism
$\mathfrak{h}:B_n\lto \R{n}$ and diffeomorphisms $\mathfrak{t}_y$, $y\in \Rm^n$,
of $\Rm^n$, equal to the identity outside of the open unit ball 
$B_n$, 
such that the map $(x,y)\lmto \mathfrak{t}_y(x)$ is smooth and such that
$$
\mathfrak{h}\circ \mathfrak{t}_y=y+\mathfrak{h}
$$
 on $B_n$. This implies that $\mathfrak{t}$ is an action of the group
$\Rm^n$ on $\Rm^n$, which means that 
$\mathfrak{t}_y\circ \mathfrak{t}_{y'}=\mathfrak{t}_{y+y'}$
for each $y, y'$.
Since $\mathfrak{t}$ is smooth,  $\mathfrak{t}_0=Id$, and 
 $\mathfrak{t}_y=Id$ outside of the unit ball, the maps 
$\mathfrak{t}_y$ converge uniformly to the identity as $y\lto 0$,
and all their derivatives converge uniformly to the derivatives of the identity. 

Let us give some details on the construction of $\mathfrak{h}$
and $\mathfrak{t}$.
We set
$$
\mathfrak{h}(x)= \frac{h(\|x\|)}{\|x\|}  x,
$$
where 
$h : {[0,1[ }\to \R{}_+$ is a smooth, strictly  increasing ($h'>0$) function such that
$$\left\{
\begin{array}{lc}
h(r) = r, & 0\leqslant r \leqslant 1/3 ,\\
h(r)= \exp \left( (r-1)^{-2} \right), & 2/3 \leqslant r < 1.
\end{array}
\right.
$$
We then define $\mathfrak{t}_y$, for each $y\in \Rm^n$ by
$$
\left\{
\begin{array}{lc}
\mathfrak{t}_y(x) = \mathfrak{h}^{-1} \big(\mathfrak{h}(x)+y\big) &\; \; \text{if} \; \; x\in B_n, \\
\mathfrak{t}_y(x) =x  &\; \; \text{if} \; \; x\in \R{n}\setminus B_n.
\end{array}
\right.
$$
It is clear from these formula that 
$ \mathfrak{t}_{y+y'}=\mathfrak{t}_y\circ\mathfrak{t}_{y'}$.
The only issue is the smoothness of $\mathfrak{t}$.
Differentiating the previous group property with respect to $y'$ and taking $y'=0$ yields the following relation:
$$
\frac{\partial }{\partial  y} \mathfrak{t}_{y}=\frac{\partial }{\partial  y} \mathfrak{t}_{0} \circ \mathfrak{t}_{y}.
$$
This implies that 
$$
\mathfrak{t}_{y}(x) =x+\int_0^1 \frac{\d}{\d t} \mathfrak{t}_{ty}(x)\d t
=x+
\int_0^1  \Big ( \frac{\partial }{\partial  y} \mathfrak{t}_{ty}(x)\Big ) 
y\d t 
= x+
\int_0^1  \Big ( \frac{\partial }{\partial  y}
 \mathfrak{t}_{0}\big(\mathfrak t_{ty}(x)\big)\Big ) y\d t.
$$
In other words, the map $\mathfrak t_y$
is the time-one flow of the vector field $X_y(x):= M(x)y$,
where $M(x)=\partial_y\mathfrak{t}_y(x)_{|y=0}$.
In order to prove that the map $\mathfrak{t}$ is smooth,
it is enough to observe that the matrix
$M(x)$ depends smoothly on $x$.
This matrix can be computed,
 recalling that the gradient of the norm 
$x\mapsto \| x \|$ is $r_x:= x/\| x \|$:
$$
M(x) = \d_{\mathfrak{h}(x)} \mathfrak{h}^{-1} 
=  \frac{1}{h'(\|x\|)} r_x^{\ t} r_x + \frac{\|x \|} {h(\|x\|)}(I_n- r_x^{\ t} r_x). 
 $$
Since $1/h$, $1/h'$, as well as all 
their derivatives go to $0$ when $\| x \| \to 1$, we conclude that 
$M(x)$ is smooth.

We have exposed the construction of  $\mathfrak{h}$ and $\mathfrak{t}$.
They allow to define a local regularization procedure 
with the help of a smooth kernel
 $K_1:\R{n}\to [0,\infty)$.
We assume that $K_1$ is supported 
in the unit ball $B_n$,  and that $\int K_1 =1$. 
For  $\eta >0$, we set  $K_{\eta} (x)=\eta^{-n} K_1 ( \eta ^{-1} x)$.

\begin{lm}\label{boule}
Let $O\subset \R{n}$ be an open set containing $\overline{B}_n$.
Given a locally integrable function $f:O\lto \Rm$ and $\eta\in ]0,1[$, we define
$$
f_\eta (x) =\int_{\R{n}} f\big(\mathfrak{t}_y(x)\big)K_\eta (-y) \d y.
$$
The following assertions hold:
\begin{enumerate}
\item The function $f_\eta$ is $C^{\infty}$ in $B_n$, and equal to $f$ outside of $B_n$,
\item If $f$ is $C^k$ on $O$, then so are the functions $f_{\eta}$, and $f_{\eta}\lto f$
in $C^k$ as $\eta\lto 0$.
\item If $f$ is $C^{k,1}$ on $O$, then so are the functions $f_{\eta}$, and
 $\limsup_{\eta\lto 0} {\rm Lip}(\d^kf_{\eta})\leq {\rm Lip}(\d^kf)$.
\item If, in some open set $O'\subset O$,  $f$ is $C^l$ in $O'$, then so is $f_\eta$.
\end{enumerate}
\end{lm}

\proof
On $B_n$ we have
$$
f_{\eta}\circ \mathfrak{h}^{-1}=(f\circ \mathfrak{h}^{-1})\star K_{\eta},
$$
where $\star$ is the convolution.
Since the functions $K_{\eta}$ are smooth, this  implies the first  claim.
Writing
$$
f_\eta -f =\int_{B(0,\eta)} (f\circ \mathfrak{t}_y-f)K_\eta (-y) \d y
$$
and observing that $f\circ \mathfrak t_y-f\lto 0$ in $C^k(\R{n},\R{n})$ as $y\lto 0$
\big(because $ \mathfrak{t}_y\lto Id$ in $C^k(\R{n},\R{n})$\big) yields the second claim.
We will now prove that 
\begin{equation}\label{limsup}
\limsup_{y\lto 0} \text{Lip}\big(\d^k(f\circ \mathfrak{t}_y)\big)\leq \text{Lip}(\d^kf),
\end{equation}
which yields the third claim in view of the relation
$$
\d^k_xf_\eta  =\int_{B(0,\eta)} \d_x^k(f\circ \mathfrak{t}_y)K_\eta (-y) \d y.
$$
Let us consider a component $\partial^{\alpha}_x(f\circ  \mathfrak{t}_y)$
of the  differential $\d^k(f\circ \mathfrak{t}_y)$, where 
$\alpha = (\alpha_1,\dots ,\alpha_n)$ is a multi--index such that
 $|\alpha|=\sum \alpha_i =k$.
By the Fa\`a di Bruno formula,
expressed in terms of partial differentials
(see \cite{faamulti} for example), we have 
$$
 \partial^\alpha_x (f\circ \mathfrak{t}_y) 
= \sum_{1\leqslant |\lambda|\leqslant  |\alpha|}\partial_{\mathfrak{t}_y(x)}^\lambda f
\cdot B_{\alpha,\lambda}(\d_x \mathfrak{t}_y,\ldots , \d_x^{|\alpha|}\mathfrak{t}_y),
$$
where the $B_{\alpha,\lambda}$ are universal multi--variable polynomials
  with no constant terms.
These polynomials satisfy the equalities
$$
B_{\alpha,\alpha}(Id,0,\cdots,0)=1\quad\text{and}\quad
B_{\alpha,\lambda }(Id,0,\cdots,0)=0
$$
for all $\lambda \neq \alpha$. 
Since $\mathfrak{t}_y\lto Id$ in $C^{\infty}$, the first of these equalities
implies that the function 
$x\lmto B_{\alpha,\alpha}(\d_x \mathfrak{t}_y,\ldots , \d_x^{|\alpha|}\mathfrak{t}_y)$
is converging to $1$ in $C^{\infty}$. 
Concerning the other factor in this term, we have 
$$
\text{Lip} \big((\partial ^\alpha f)\circ \mathfrak{t}_y\big)
\leq \text{Lip}(\partial^{\alpha }f) \text{Lip}(\mathfrak{t}_y)
\lto \text{Lip}(\partial^{\alpha }f).
$$
We deduce that the upper limit of the Lipschitz constants of the term
corresponding to $\lambda=\alpha$ is not greater than $\text{Lip}(\partial^{\alpha }f)$.

On the other hand, for each of the terms with $\lambda \neq \alpha$,
the function 
$x\lmto B_{\alpha,\lambda}(\d_x \mathfrak{t}_y,\ldots , \d_x^{|\alpha|}\mathfrak{t}_y)
$
is converging to $0$ in $C^{\infty}$ hence the Lipschitz constant of the 
function 
$$
x\lmto \partial_{\mathfrak{t}_y(x)}^\lambda f
\cdot B_{\alpha,\lambda}(\d_x \mathfrak{t}_y,\ldots , \d_x^{|\alpha|}\mathfrak{t}_y)
$$
is converging to $0$. We conclude that 
$$
\limsup \text{Lip} \big(\partial^{\alpha}(f\circ \mathfrak{t}_y )\big)
\leq \text{Lip} (\partial^{\alpha}f),
$$
which implies (\ref{limsup}) hence the third point of the statement.

Regarding the last claim of the statement, we 
consider the set
$
\Omega:=\cap_{y\in \overline B(0,\eta)}
\mathfrak{t}_y ^{-1}(O'),
$
and claim that $\Omega$ is open.
Assuming the claim, we observe that the function 
$f_{\eta}$ is smooth in $B_n$ and that it is $C^l$ 
in $\Omega$.
Since the maps $\mathfrak{t}_y$ are all the identity outside of $B_n$,
the set $\Omega$ contains $O'-B_n$.
We have covered $O'$ by two open sets, $B_n$ and $\Omega$, such that the $f_{\eta}$
is $C^l$ on each of them, we conclude that this function
is $C^l$ on $O'$.

To prove that $\Omega$ is open,
we fix $x_0\in \Omega$. 
For each $y_0\in \overline B(0,\eta)$, we have 
$\mathfrak{t}_{y_0} (x_0)\in  O'$,
hence there exists an open set $U_{y_0}$ containing $y_0$  and an open set 
$\Omega_{y_0}$ 
containing $x_0$ such that
$\mathfrak{t}_y (x)\in  O'$ far all $(x,y)\in \Omega_{y_0}\times U_{y_0}$.
By compactness, there exists finitely many points $y_i\in \overline B(0,\eta)$
such that the open sets $U_{y_i}$ cover $\overline B(0,\eta)$.
The open intersection $\cap \Omega_{y_i}$,
which contains $x_0$, is then contained in $\Omega$.
Since this holds for each $x_0\in \Omega$, we have proved that $\Omega$ is open.
 \qed

\begin{lm}\label{ouvert}
Let $O$ be open subsets of $\R{n}$
and let $f: O\to \R{}$ be a $C^{k,1}$ function.
Given a continuous function  $\e : O \to [0,\infty)$,
there exists a function $f_\e$ such that:
 \begin{enumerate}
\item the function $f_\e$ is $C^{\infty}$ in the open set 
 $\{x\in O ,\e(x) >0\}\subset O$,
\item    
$|f_{\e}(x)-f(x)|+\| \d_x f_{\e} -\d_x f\|+\cdots+\| \d^k_x f_{\e} -\d^k_x f\|
\leqslant \e(x)$ for each  $x\in O$,
\item  the function $f_{\e}$
is $C^{k,1}$ on $O$, and
${\rm Lip}(\d^kf_{\e})\leq 1+{\rm Lip}(\d^kf).
$
\end{enumerate}

\end{lm}
\proof
Let us denote by $F$ the closed set $\{\e=0\}$.
The complement of $F$ in $O$ is open, and 
we consider a  locally finite covering $(O_i)_{i\in \mathbb{N}^*}$ of $O \setminus F$ 
by  open balls  compactly included in $O \setminus F$. 
Since  $\inf \{ \e(x),\quad x\in O_i\}>0$. we can 
 construct inductively, using  Lemma \ref{boule} a sequence of functions, 
$(f_i)_{i\in \mathbb{N}}$ such that 
\begin{itemize}
\item $f_0=f$,
\item for each $i\in \mathbb{N}$, the function $f_{i+1}$ is $C^{\infty}$ in 
$O_1\cup \cdots \cup O_{i+1}$,
\item for each $i\in \mathbb{N}$, 
 the functions $f_i$ and $f_{i+1}$ are equal in $O \setminus O_{i+1}$,
\item for each $i\in \mathbb{N}$, the function $f_{i+1}$ is $C^{k,1}$ in $O$,
and ${\rm Lip}(\d^kf_{i+1})\leq 2^{-i-1}+{\rm Lip}(\d^kf_i),
$
\item 
$|f_{i+1}(x)-f_i(x)|+\| \d_x f_{i+1} -\d_x f_i\|+\cdots+\| \d^k_x f_{i+1} -\d^k_x f_i\|
\leqslant 2 ^{-1-i}\e(x)$ for each  $x\in O$, $ i\in \Nm$,
\end{itemize}
Each point of $O$ has a neighborhood on which the sequence $f_i$
is eventually constant, hence the limit $f_{\e}:= \lim f_i$
is well-defined and smooth  on $\cup_i O_i=O\setminus F$.
The desired estimates on $f_{\e}$ follow immediately
from the inductive estimates by summation.  
\qed

\textsc{Proof of Theorem \ref{approx}}.
We fix a locally finite atlas $(\phi_i)_{i\in \mathbb{N}^*}$ constituted of 
smooth maps $\phi_i : 2B_n\to M$, where $B_n$ is the open unit ball. We assume  
that all the images $\phi_i(2B_n)$, $i\in \mathbb{N}^*$ are relatively compact in $M$ and that the $\phi_i(B_n)$, $i\in \mathbb{N}^*$ still cover $M$.
By Lemma \ref{ouvert}, it is possible to construct inductively a sequence of functions $f_i$, by iteratively modifying $f_i\circ \phi_{i+1}$ on $B_n$, such that 
\begin{itemize}
\item $f_0=f$,
\item for each $i\in \mathbb{N}$, the function $f_{i+1}$ is $C^{\infty}$ 
in $\bigcup_{j\leqslant i+1}\phi_{j}(B_n) \cap \Omega$,
\item for each $i\in \mathbb{N}$, in $M \setminus \phi_{i+1}(B_n)$, the functions $f_i$ and $f_{i+1}$ are equal,
\item for each $i\in \mathbb{N}$, the function $f_{i+1}$ is $C^{k,1}$ on $M$,
\item for each $ i\in \mathbb{N}$, $ x\in M $,
$ |f_i(x)-f_{i+1}(x)|+
\cdots +\|\d^k_{x} f_i  - \d^k_{x} f_{i+1}\|\leqslant  2^{-i-1}\e(x).$
\end{itemize}

Each point $x\in M$ has a neighborhood on which the sequence $f_i$
is eventually constant, hence the limit $g=\lim f_i$
is well defined, locally $C^{k,1}$, and smooth on $\Omega$.
The inequality on the differentials follows by summation from the iterative 
assumptions.
\qed

\bibliography{CC-09-08}

\begin{thebibliography}{CISM00}

\bibitem[AD00]{Amb}
L.~Ambrosio and N.~Dancer.
\newblock {\em Calculus of variations and partial differential equations}.
\newblock Springer-Verlag, Berlin, 2000.
\newblock Topics on geometrical evolution problems and degree theory, Papers
  from the Summer School held in Pisa, September 1996, Edited by G. Buttazzo,
  A. Marino and M. K. V. Murthy.

\bibitem[BB07]{Be}
P.~Bernard and B.~Buffoni.
\newblock Weak {KAM} pairs and {M}onge-{K}antorovich duality.
\newblock In {\em Asymptotic analysis and singularities---elliptic and
  parabolic {PDE}s and related problems}, volume~47 of {\em Adv. Stud. Pure
  Math.}, pages 397--420. Math. Soc. Japan, Tokyo, 2007.

\bibitem[Ber02]{bernard2002connecting}
P.~Bernard.
\newblock Connecting orbits of time dependent {L}agrangian systems.
\newblock {\em Annales de l'institut Fourier}, 52(5):1533--1568, 2002.

\bibitem[Ber07]{Be1}
P.~Bernard.
\newblock Existence of {$C\sp {1,1}$} critical sub-solutions of the
  {H}amilton-{J}acobi equation on compact manifolds.
\newblock {\em Ann. Sci. \'Ecole Norm. Sup. (4)}, 40(3):445--452, 2007.

\bibitem[Ber08]{Be111}
P.~Bernard.
\newblock The dynamics of pseudographs in convex {H}amiltonian systems.
\newblock {\em J. Amer. Math. Soc.}, 21(3):615--669, 2008.

\bibitem[Ber10]{Beil}
P.~Bernard.
\newblock Lasry-{L}ions regularization and a lemma of {I}lmanen.
\newblock {\em Rend. Semin. Mat. Univ. Padova}, 124:221--229, 2010.

\bibitem[BR05]{bernard2005convergence}
P.~Bernard and J.M. Roquejoffre.
\newblock Convergence to time-periodic solutions in time-periodic
  {H}amilton--{J}acobi equations on the circle.
\newblock {\em Communications in Partial Differential Equations},
  29(3-4):457--469, 2005.

\bibitem[Car01]{Carda}
P.~Cardaliaguet.
\newblock Front propagation problems with nonlocal terms. {II}.
\newblock {\em J. Math. Anal. Appl.}, 260(2):572--601, 2001.

\bibitem[CISM00]{cis}
G.~Contreras, R.~Iturriaga, and H.~Sanchez-Morgado.
\newblock {W}eak solutions of the {H}amilton-{J}acobi equation for time
  periodic {L}agrangians.
\newblock {\em preprint}, 2000.

\bibitem[CS96]{faamulti}
G.~M. Constantine and T.~H. Savits.
\newblock A multivariate {F}a\`a di {B}runo formula with applications.
\newblock {\em Trans. Amer. Math. Soc.}, 348(2):503--520, 1996.

\bibitem[dR73]{derh}
G.~de~Rham.
\newblock {\em Vari\'et\'es diff\'erentiables. {F}ormes, courants, formes
  harmoniques}.
\newblock Hermann, Paris, 1973.
\newblock Troisi{\`e}me {\'e}dition revue et augment{\'e}e, Publications de
  l'Institut de Math{\'e}matique de l'Universit{\'e} de Nancago, III,
  Actualit{\'e}s Scientifiques et Industrielles, No. 1222b.

\bibitem[Fat98]{fathi1998convergence}
A.~Fathi.
\newblock Sur la convergence du semi-groupe de {L}ax-{O}leinik.
\newblock {\em Comptes Rendus de l'Acad{\'e}mie des Sciences-Series
  I-Mathematics}, 327(3):267--270, 1998.

\bibitem[FM00]{fathi2000failure}
A.~Fathi and J.~Mather.
\newblock Failure of convergence of the {L}ax-{O}leinik semi-group in the
  time-periodic case.
\newblock {\em Bull. Soc. Math. France}, 128(3):473--483, 2000.

\bibitem[FS04]{FaSi}
A.~Fathi and A.~Siconolfi.
\newblock Existence of {$C\sp 1$} critical subsolutions of the
  {H}amilton-{J}acobi equation.
\newblock {\em Invent. Math.}, 155(2):363--388, 2004.

\bibitem[FZ10]{fz}
A.~Fathi and M.~Zavidovique.
\newblock Ilmanen's lemma on insertion of {$C^{1,1}$} functions.
\newblock {\em Rend. Semin. Mat. Univ. Padova}, 124:203--219, 2010.

\bibitem[Gom05]{Gom}
D.A. Gomes.
\newblock Viscosity solution method and the discrete {A}ubry-{M}ather problem.
\newblock {\em Discrete and Continuous Dynamical Systems, Series A},
  13:103--116, 2005.

\bibitem[Hir94]{Hir}
M.~W. Hirsch.
\newblock {\em Differential topology}, volume~33 of {\em Graduate Texts in
  Mathematics}.
\newblock Springer-Verlag, New York, 1994.
\newblock Corrected reprint of the 1976 original.

\bibitem[Ilm93]{Il}
T.~Ilmanen.
\newblock The level-set flow on a manifold.
\newblock In {\em Differential geometry: partial differential equations on
  manifolds ({L}os {A}ngeles, {CA}, 1990)}, volume~54 of {\em Proc. Sympos.
  Pure Math.}, pages 193--204. Amer. Math. Soc., Providence, RI, 1993.

\bibitem[Zav10]{Za2}
M.~Zavidovique.
\newblock Existence of {$C^{1,1}$} critical subsolutions in discrete weak {KAM}
  theory.
\newblock {\em J. Mod. Dyn.}, 4(4):693--714, 2010.

\bibitem[Zav12]{Za}
M.~Zavidovique.
\newblock Strict subsolutions and {M}a\~n\' e potential in discrete weak {KAM}
  theory.
\newblock {\em {C}ommentarii {M}athematici {H}elvetici}, 87(1):1--39, 2012.

\end{thebibliography}
\bibliographystyle{alpha}

\end{document}